\documentstyle[twoside,11pt]{article}
\date{}

\title{Principal Values and Principal Subspaces of Two Subspaces \\ of Vector Spaces with Inner Product}

\author{Ice B. Risteski$^1$, Kostadin G. Tren\v{c}evski$^2$\\
$^1$2 Milepost Place \#606, Toronto, ON, M4H 1C7 Canada,\\ e-mail: iceristeski@hotmail.com
\\
$^2$Institute of Mathematics, Sts. Cyril and Methodius University in Skopje,  
\\ 
P.O.Box 162, 1000 Skopje, Macedonia, e-mail: kostatre@pmf.ukim.mk}

\begin{document}
\maketitle

\begin{abstract}
In this paper
is studied the problem concerning the angle between two 
subspaces of arbitrary dimensions in Euclidean space $E_{n}$. 
It is proven that the angle between two subspaces is equal to the 
angle between their orthogonal subspaces. Using 
the eigenvalues and eigenvectors of corresponding matrix representations, 
there are introduced principal values and principal subspaces. Their 
geometrical interpretation is also given together with the canonical 
representation of the two subspaces. The canonical matrix for the 
two subspaces is introduced and its properties of duality are obtained. 
Here obtained results expand the classic results given in 
[1,2]. 


\end{abstract}
\vskip 1cm

\noindent {\bf 1. Angle between two subspaces in ${\bf E_{n}}$}
\medskip

\noindent 
We will prove the following theorem which will enable us to define angle 
between two subspaces of arbitrary dimensions of the Euclidean space 
$E_{n}$.

\noindent{\bf Theorem 1.1.} {\em Let ${\bf a}_{1},\cdots ,{\bf a}_{p}$ and 
${\bf b}_{1},\cdots ,{\bf b}_{q}$ are bases of two subspaces 
$\Sigma _{1}$ and $\Sigma _{2}$ of Euclidean space
$E_{n}$ with inner product $( , )$ respectively 
and suppose that $p \le q \le n$. 
If $p<q$, assume that ${\bf b}_{1},\cdots ,{\bf b}_{q}$ is an orthonormal
base. Then the following inequality holds} 
$$det [MM^{T}]\le 
\left \vert \matrix{ 
({\bf a}_{1}, {\bf a}_{1})&({\bf a}_{1}, {\bf a}_{2})& \cdots &
({\bf a}_{1}, {\bf a}_{p})\cr 
({\bf a}_{2}, {\bf a}_{1})&({\bf a}_{2}, {\bf a}_{2})& \cdots &
({\bf a}_{2}, {\bf a}_{p})\cr 
\cdot & & & \cr 
\cdot & & & \cr 
\cdot & & & \cr 
({\bf a}_{p}, {\bf a}_{1})&({\bf a}_{p}, {\bf a}_{2})& \cdots &
({\bf a}_{p}, {\bf a}_{p})\cr }\right \vert \times \leqno{(1.1)}$$
$$\times \left \vert \matrix{ 
({\bf b}_{1}, {\bf b}_{1})&({\bf b}_{1}, {\bf b}_{2})& \cdots &
({\bf b}_{1}, {\bf b}_{q})\cr 
({\bf b}_{2}, {\bf b}_{1})&({\bf b}_{2}, {\bf b}_{2})& \cdots &
({\bf b}_{2}, {\bf b}_{q})\cr 
\cdot & & & \cr 
\cdot & & & \cr 
\cdot & & & \cr 
({\bf b}_{q}, {\bf b}_{1})&({\bf b}_{q}, {\bf b}_{2})& \cdots &
({\bf b}_{q}, {\bf b}_{q})\cr }\right \vert ,$$
{\em where }
$$ M=\left [ \matrix{ 
({\bf a}_{1}, {\bf b}_{1})&({\bf a}_{1}, {\bf b}_{2})& \cdots &
({\bf a}_{1}, {\bf b}_{q})\cr 
({\bf a}_{2}, {\bf b}_{1})&({\bf a}_{2}, {\bf b}_{2})& \cdots &
({\bf a}_{2}, {\bf b}_{q})\cr 
\cdot & & & \cr 
\cdot & & & \cr 
\cdot & & & \cr 
({\bf a}_{p}, {\bf b}_{1})&({\bf a}_{p}, {\bf b}_{2})& \cdots &
({\bf a}_{p}, {\bf b}_{q})\cr }\right ] $$
{\em and moreover equality holds if and only if 
$\Sigma _{1}$ is subspace of $\Sigma _{2}$.}

\medskip
\noindent {\em Proof.} 
The inequality (1.1) is invariant under any elementary 
row operation of ${\bf a}_{1},\cdots ,{\bf a}_{p}$, and if $p=q$ 
for any elementary row operation of ${\bf b}_{1},\cdots ,{\bf b}_{q}$. 
Thus, without loss of generality we can assume that 
$\{ {\bf a}_{1}, \cdots ,{\bf a}_{p}\}$ and 
$\{ {\bf b}_{1}, \cdots ,{\bf b}_{q}\}$ are orthonormal systems. 
Then we should prove that 
$$ det [MM^{T}] \le 1.$$
Let denote 
$$ {\bf c}_{i} = (({\bf a}_{i}, {\bf b}_{1}),
({\bf a}_{i}, {\bf b}_{2}),\cdots ,({\bf a}_{i}, {\bf b}_{q}))
\in {\bf R}^{q}\quad (1\le i\le p).$$
Since $\{ {\bf b}_{i}\}$ and $\{ {\bf a}_{i}\}$ are orthonormal systems 
we get that $\Vert {\bf c}_{i}\Vert \le 1$ with respect to the 
Euclidean metric in ${\bf R}^{q}$. 

Let ${\bf c}_{p+1},
\cdots ,{\bf c}_{q}$ be an orthonormal system of vectors
such that each of them is orthogonal to 
${\bf c}_{1},\cdots ,{\bf c}_{p}$. Then 
$$ det [MM^{T}] = 
\left \vert \matrix{ 
({\bf c}_{1}\cdot {\bf c}_{1})&({\bf c}_{1}\cdot {\bf c}_{2})&\cdots &
({\bf c}_{1}\cdot {\bf c}_{p})\cr 
({\bf c}_{2}\cdot {\bf c}_{1})&({\bf c}_{2}\cdot {\bf c}_{2})&\cdots &
({\bf c}_{2}\cdot {\bf c}_{p})\cr 
\cdot & & & \cr
\cdot & & & \cr
\cdot & & & \cr
({\bf c}_{p}\cdot {\bf c}_{1})&({\bf c}_{p}\cdot {\bf c}_{2})&\cdots&
({\bf c}_{p}\cdot {\bf c}_{p})\cr }\right \vert =$$
$$= \left \vert \matrix{ 
({\bf c}_{1}\cdot {\bf c}_{1})&({\bf c}_{1}\cdot {\bf c}_{2})&\cdots &
({\bf c}_{1}\cdot {\bf c}_{q})\cr 
({\bf c}_{2}\cdot {\bf c}_{1})&({\bf c}_{2}\cdot {\bf c}_{2})&\cdots &
({\bf c}_{2}\cdot {\bf c}_{q})\cr 
\cdot & & & \cr
\cdot & & & \cr
\cdot & & & \cr
({\bf c}_{q}\cdot {\bf c}_{1})&({\bf c}_{q}\cdot {\bf c}_{2})&\cdots&
({\bf c}_{q}\cdot {\bf c}_{q})\cr }\right \vert $$
which is the square of the volume of the parallelotop in 
${\bf R}^{q}$ generated by the 
vectors ${\bf c}_{1},\cdots ,{\bf c}_{q}$. Since 
$\Vert {\bf c}_{i}\Vert \le 1$, $(1\le i\le q)$ we obtain 
$det [MM^{T}]\le 1$. 

Moreover, equality holds if and only if 
${\bf c}_{1},\cdots ,{\bf c}_{q}$ 
is an orthonormal system. But $\Vert {\bf c}_{i}\Vert =1$ implies 
that ${\bf a}_{i}$ belongs to the subspace $\Sigma _{2}$. Thus 
$\Sigma _{1}\subseteq \Sigma _{2}$. Conversely, if 
$\Sigma _{1}\subseteq \Sigma _{2}$ then it is trivial that equality 
holds in (1.1). \quad {\hspace* {\fill}{$\Box $}}
\medskip

\noindent Under the assumptions of Theorem 1.1 we define 
the angle $\varphi$
between $\Sigma _{1}$ and $\Sigma _{2}$ by 
$$ \cos \varphi = {\sqrt{det [MM^{T}]}\over 
\sqrt{\Gamma _{1}}\cdot \sqrt{\Gamma _{2}}}\leqno{(1.2)}$$
where the matrix $M$ was defined in Theorem 1.1 and $\Gamma _{1}$ and 
$\Gamma _{2}$ are the Gram's determinants obtained by the vectors 
${\bf a}_{1},\cdots ,{\bf a}_{p}$ and 
${\bf b}_{1},\cdots ,{\bf b}_{q}$ respectively. 

Note that $det [MM^{T}]\ge 0$; considering both values of 
$\sqrt{det [MM^{T}]}$, we obtain two angles $\varphi $ and $\pi -\varphi $. 
Note that $det [MM^{T}]=0$ if $q <p $. 

In this paper we give some deeper results concerning the Theorem 1.1.
Indeed, some theorems which yield to 
principal directions on both subspaces $\Sigma _{1}$ and 
$\Sigma _{2}$ and common principal values are proven. 

\medskip
\noindent In the next research will be used the following result. 

\medskip
\noindent 
{\bf Theorem 1.2.} {\em Let $U$ be any $p \times q $ matrix. 
Any nonzero scalar $\lambda $ is an eigenvalue of the square matrix 
$UU^{T}$ if and only if it is eigenvalue of the square matrix 
$U^{T}U$ and moreover the multiplicities of $\lambda $ for both 
matrices $UU^{T}$ and $U^{T}U$ are equal. }

\medskip
\noindent {\em Proof.} 
Assume that $\lambda \neq 0$ is an eigenvalue of $UU^{T}$ 
with geometrical multiplicity $r$ and assume that 
${\bf x}_{1},\cdots ,{\bf x}_{r}$ are linearly independent 
eigenvectors corresponding to $\lambda $. Then we will prove that the vectors
$${\bf y}_{i}=U^{T}{\bf x}_{i} ,\quad (1\le i\le r)$$
are linearly independent eigenvectors for the matrix $U^{T}U$. Indeed, 
$$U^{T}U{\bf y}_{i}=(U^{T}U)U^{T}{\bf x}_{i}
=U^{T}(UU^{T}{\bf x}_{i})=
\lambda U^{T}{\bf x}_{i}= \lambda {\bf y}_{i}$$
and thus ${\bf y}_{i}$ are eigenvectors of $U^{T}U$ 
corresponding to the eigenvalue $\lambda $. 

Now let us assume that 
$\alpha _{1}{\bf y}_{1}+\cdots +\alpha _{r}{\bf y}_{r}=0,$ then 
multiplying this equality by $U$ from left we obtain
$$\lambda \alpha _{1}{\bf x}_{1}+\cdots +\lambda \alpha _{r}{\bf x}_{r}=0.$$
Since $\lambda \neq 0$ we obtain 
$$\alpha _{1}{\bf x}_{1}+\cdots +\alpha _{r}{\bf x}_{r}=0$$
and hence $\alpha _{1}=\cdots =\alpha _{r}=0$ because 
${\bf x}_{1},\cdots ,{\bf x}_{r}$ are linearly independent vectors. 

Hence the geometric multiplicity of $\lambda $ for the matrix 
$UU^{T}$ is smaller or equal to the 
geometric multiplicity of $\lambda $ for the matrix 
$U^{T}U$. Analogously, 
the geometric multiplicity of $\lambda $ for the matrix 
$U^{T}U$ is smaller or equal to the 
geometric multiplicity of $\lambda $ for the matrix $UU^{T}$. 
Thus these two geometrical multiplicities are equal. 
Since $UU^{T}$ and $U^{T}U$ are symmetric 
non-negative definite matrices, we obtain that their geometrical 
multiplicities are equal to the algebraic multiplicities. 
\quad {\hspace* {\fill}{$\Box $}}
\medskip

\noindent 
Now we are enable to prove the following theorem. 

\medskip
\noindent 
{\bf Theorem 1.3.} {\em If $\Sigma _1 $ and $\Sigma _2 $ are any subspaces
of the Euclidean vector space $E_{n}$ and 
$\Sigma _1 ^* $ and $\Sigma _2 ^*$ are their orthogonal complements, then} 
$$\varphi(\Sigma _1 ,\Sigma _2 )=\varphi(\Sigma _1^* ,\Sigma _2^* ).$$

\medskip
\noindent {\em Proof.} 
Assume that $dim \Sigma _{1}=p$ and $dim \Sigma _{2}=q$. 
Without loss of generality we assume that $p \le q $ and assume that 
$\Sigma _1 $ is generated by ${\bf e}_i,$ $(1\le i\le p )$ and 
$\Sigma _1^* $ is generated by ${\bf e}_j,$ $(p +1\le j\le n)$ 
where ${\bf e}_i,$ $(1\le i\le n)$ is the standard basis of $E_{n}$. 
Further without loss of generality we can assume that 
$\Sigma _2 $ is generated by ${\bf a}_i$, $(1\le i\le q )$ and 
$\Sigma _2^* $ is generated by ${\bf a}_j$, $(q +1\le j\le n)$, 
where ${\bf a}_i$, $(1\le i\le n)$ is an orthonormal system of vectors. 
Let ${\bf a}_i $ has coordinates 
$(a_{i1},a_{i2},\cdots ,a_{in})$, $(1\le i\le n)$ and the matrix 
with row vectors ${\bf a}_{1},\cdots ,{\bf a}_{n}$ will be denoted by $A$. 
We denote by $X$, $Y$ and $Z$ the following submatrices of A: 
$X$ is submatrix of $A$ with elements $a_{ij}$, $(1\le i\le p;\; 
1\le j\le q)$; 
$Y$ is submatrix of $A$ with elements $a_{ij}$, $(1\le i\le p;\; 
q+1\le j\le n)$; 
$Z$ is the submatrix of $A$ with elements $a_{ij}$, $(p +1\le i\le n;\; 
q +1\le j\le n)$. According to these assumptions 
$$\cos ^{2}\varphi(\Sigma _1 ,\Sigma _2 ) = det [XX^{T}]$$ 
and 
$$\cos ^{2}\varphi(\Sigma _1^* ,\Sigma _2^* ) = det [Z^{T}Z]$$
and we should prove that 
$$ det [XX^{T}] = det [Z^{T}Z] .$$
Since $A$ is an orthogonal matrix, it holds 
$$XX^{T} = I_{p \times p } - YY^{T}\quad \hbox { and } \quad 
Z^{T}Z = I_{(n-q) \times (n-q)} - Y^{T}Y$$
and we should prove that 
$$det [I_{p \times p } - YY^{T}] = 
det [I_{(n-q) \times (n-q)} - Y^{T}Y].$$
Let $\lambda _{1},\cdots ,\lambda _{p}$ be the eigenvalues of $YY^{T}$ 
and $\mu _{1},\cdots ,\mu _{n-q}$ be the eigenvalues of $Y^{T}Y$. 
According to Theorem 1.2, the matrices $YY^{T}$ and $Y^{T}Y$ have 
the same non-zero eigenvalues with the same multiplicities and hence
$$ det [I_{p \times p } - YY^{T}] = 
(1-\lambda _{1})\cdots (1-\lambda _{p}) =$$
$$ = (1-\mu _{1})\cdots (1-\mu _{q}) = 
det [I_{(n-q) \times (n-q)} - Y^{T}Y]. \eqno{\Box}$$

\medskip
\noindent {\bf 2. Principal values and principal subspaces }
\medskip

\noindent First we prove the following statement. 

\medskip
\noindent 
{\bf Theorem 2.1.} {\em Let $\Sigma _{1}$ and $\Sigma _{2}$ be two 
vector subspaces of the Euclidean space $E_{n}$ of dimensions 
$p$ and $q$, $(p\le q)$ and let $A_{1}$ and $A_2$ be 
$p\times n$ and $q\times n$ 
matrices whose vector rows generate the subspace $\Sigma _1$ 
and $\Sigma _2$ respectively. Then the eigenvalues of the matrix }
$$ f(A_{1},A_{2}) = A_{1}A_{2}^{T}(A_{2}A_{2}^{T})^{-1}
A_{2}A_{1}^{T}(A_{1}A_{1}^{T})^{-1}$$
{\em are $p$ canonical squares 
$\cos ^{2}\varphi _{i}$, $(1\le i\le p )$ and moreover }
$$\cos ^{2}\varphi = \prod _{i=1}^{p} \cos ^{2}\varphi _{i},$$
{\em where $\varphi $ is the angle between the subspaces $\Sigma _{1}$ and }
$\Sigma _{2}$.

\medskip
\noindent {\em Proof.} 
The transition of the base of $\Sigma _{j}$ to another base 
corresponds to multiplication of $A_{j}$ by nonsingular matrix $P_{j}$, 
i.e. $A_{j}\rightarrow P_{j}A_{j}$, where $P_1$ is $p\times p$ matrix and 
$P_{2}$ is $q\times q$ matrix. By direct calculation one verifies that 
$$ f(P_{1}A_{1}, P_{2}A_{2}) = P_{1} f(A_{1},A_{2}) P_{1}^{-1}$$
and thus the eigenvalues are unchanged. Moreover, $f(A_{1},A_{2})$ 
are unchanged under the transformation of form $A_{j}\rightarrow A_{j}R$
where $R$ is any orthogonal matrix of $n$-th order, which means that 
$f(A_{1},A_{2})$ is invariant under the change of the rectangular 
Cartesian coordinates in the Euclidean space $E_{n}$. 

Since $A_{1}A_{1}^{T}$ and $A_{2}A_{2}^{T}$ 
are positive definite matrices, there exist 
symmetric positive definite 
matrices $P_{1}$ and $P_{2}$ of orders $p$ and $q$ respectively such that 
$$P_{1}A_{1}A_{1}^{T}P_{1}^{T} = B_{1}B_{1}^{T} = I_{p\times p}
\quad \hbox { and } \quad 
P_{2}A_{2}A_{2}^{T}P_{2}^{T} = B_{2}B_{2}^{T} = I_{q\times q},$$
where $B_{1}$ and $B_2$ correspond to another bases of $\Sigma _{1}$ and 
$\Sigma _2$. Since $S=(B_{1}B_{2}^{T})(B_{1}B_{2}^{T})^{T}$
is non-negative definite matrix, there exists a 
symmetric non-negative definite orthogonal matrix $Q_{1}$ of order 
$p$ such that $Q_{1}SQ_{1}^{-1}$ is diagonalized, i.e. 
$$Q_{1}SQ_{1}^{-1} = (C_{1}B_{2}^{T})(C_{1}B_{2}^{T})^{T} = 
diag(c_{1}^{2}, c_{2}^{2}, \cdots ,c_{p}^{2}), \quad 
(c_{1}\ge c_{2}\ge \cdots \ge c_{p}\ge 0)$$
where $C_{1}=Q_{1}B_{1}$ corresponds to another basis of $\Sigma _{1}$. 
Having in mind that each $c_{i}$ is an inner product of two 
unimodular vectors, we get $c_{i}=\cos \varphi _{i}$, 
$0\le \varphi _{1}\le \varphi _{2}\le \cdots \le \varphi _{p}
\le \pi /2$. The vector rows of $C_{1}B_{2}^{T}$ are mutually orthogonal, 
which means that there exists an orthogonal matrix $Q_{2}$ of order 
$q$, such that 
$$C_{1}B_{2}^{T}Q_{2}^{T} = C_{1}C_{2}^{T} = 
\cos \varphi_{i} \delta _{ik} ,$$
where $C_{2}=Q_{2}B_{2}$ corresponds to another orthonormal base 
of $\Sigma _2 $. This shows that the ordered set of angles 
$\varphi _{1},\varphi _{2}, \cdots ,\varphi _{p}$ is canonical and 
its invariance follows from the decomposition
$$ det[\lambda I_{p\times p}-f(C_{1},C_{2})] = \prod _{i=1}^{p}
(\lambda -\cos ^{2}\varphi _{i}) = 
det[\lambda I_{p\times p}-f(A_{1},A_{2})]. $$

Finally note that according to the chosen bases of $\Sigma _1$ and 
$\Sigma _{2}$, we obtain 
$$\cos ^{2}\varphi = det [f(C_{1},C_{2})] = det [f(A_{1},A_{2})] =
\prod _{i=1}^{p} \cos ^{2}\varphi _{i}$$
where $\varphi $ is the angle between the subspaces $\Sigma _{1}$ and 
$\Sigma _{2}$. \quad {\hspace* {\fill}{$\Box $}}
\medskip

\noindent
Note that if the bases of $\Sigma _1$ and $\Sigma _2$ are orthonormal, 
then $A_{1}A_{1}^{T}=A_{2}A_{2}^{T}=I$ and 
$f(A_{1},A_{2})=A_{1}A_{2}^{T}(A_{1}A_{2}^{T})^{T}$. 

Now let us consider the case $p\ge q$. Instead of the matrix 
$f(A_1 ,A_2 )$ we should consider the matrix $f(A_2 ,A_1 )$ which is 
of type $q \times q $. Analogously to Theorem 2.1 
the eigenvalues of $f(A_2 ,A_1 )$ 
are $q$ canonical squares of cosine functions but the product 
of them is equal to zero if $p>q$. Now we prove the following 
theorem considering the mutually eigenvalues of 
$f(A_1 ,A_2 )$ and $f(A_2 ,A_1 )$. 

\medskip
\noindent
{\bf Theorem 2.2.} {\em Any nonzero scalar $\lambda $ is an eigenvalue 
of $f(A_1 ,A_2 )$ if and only if it is eigenvalue of $f(A_2 ,A_1 )$ 
and moreover the multiplicities of $\lambda $ for both 
matrices $f(A_1 ,A_2 )$ and $f(A_2 ,A_1 )$ are equal. }

\medskip
\noindent {\em Proof.} 
Let $C_{1}$ and $C_{2}$ have the same meanings like in the 
Theorem 2.1. According to Theorem 1.2 we obtain that 
any nonzero scalar $\lambda $ is an eigenvalue 
of $f(C_1 ,C_2 )$ if and only if it is eigenvalue of $f(C_2 ,C_1 )$ 
and moreover the multiplicities of $\lambda $ for both 
matrices $f(C_1 ,C_2 )$ and $f(C_2 ,C_1 )$ are equal, because 
$f(C_1 ,C_2 )=(C_{1}C_{2}^{T})(C_{1}C_{2}^{T})^{T}$. 
On the other hand, $f(A_1 ,A_2 )$ has the same eigenvalues as 
$f(C_1 ,C_2 )$ with the same multiplicities and 
$f(A_2 ,A_1 )$ has the same eigenvalues as 
$f(C_2 ,C_1 )$ with the same multiplicities. 
{\hspace* {\fill}{$\Box $}}
\medskip

\noindent
Note that $\lambda =0$ is eigenvalue for the matrix $f(A_2 ,A_1 )$ if 
$q >p $, but $\lambda =0$ may not be eigenvalue for the matrix 
$f(A_1 ,A_2 )$. 

The common eigenvalues will be called {\em principal values}. 
According to the Theorems 2.1 and 2.2 there are unique decompositions of 
the subspaces $\Sigma _{1}$ and $\Sigma _{2}$ into the orthogonal 
eigenspaces for the common non-negative eigenvalues and 
for the zero eigenvalue if such exists. These eigenspaces 
are called {\em principal subspaces} or {\em principal directions} 
for the eigenvalues with multiplicity 1. 
The geometrical interpretation of the principal values and principal 
subspaces will be given after the proof of the Theorem 2.3. 

\medskip
\noindent 
{\bf Theorem 2.3.} {\em The function 
$\cos ^{2} \varphi $, where $\varphi $ is the angle between any vector 
${\bf x}\in \Sigma _{1}$ and the subspace $\Sigma _{2}$, has maximum 
if and only if the vector ${\bf x}$ 
belongs to a principal subspace of $\Sigma _{1}$ which corresponds to 
the maximal principal value. The maximal value of 
$\cos ^{2} \varphi $ is the maximal principal value.}

\medskip
\noindent {\em Proof.} 
According to the proof of Theorem 2.1, without loss 
of generality we can suppose that $\Sigma _{1}$ is generated by the 
orthonormal vectors ${\bf a}_{i}$, $(1\le i\le p)$ and 
$\Sigma _{2}$ is generated by the 
orthonormal vectors ${\bf b}_{j}$, $(1\le j\le q)$ such that 
$({\bf a}_{i},{\bf b}_{j})=0$, $(i\neq j;\; 1\le i\le p,\; 1\le j\le q)$. 
Let ${\bf x}=\alpha _{1}{\bf a}_{1}+\cdots +\alpha _{p}{\bf a}_{p}$ 
and let $\lambda _{1}^{2}={\bf a}_{1}{\bf b}_{1}$ be the maximal 
principal value and the corresponding subspace of $\Sigma _{1}$ be 
generated by ${\bf a}_{1},\cdots ,{\bf a}_{r}$. Then for the angle 
$\varphi $ between ${\bf x}$ and $\Sigma _{2}$ it holds 
$$\cos ^{2}\varphi = 
{(\alpha _{1}\lambda _{1})^{2}+\cdots +(\alpha _{p}\lambda _{s})^{2}
\over \alpha _{1}^{2}+\cdots +\alpha _{p}^{2}} =$$
$$ = {\lambda _{1}^{2}(\alpha _{1}^{2}+\cdots +\alpha _{r}^{2})+
\lambda _{r+1}^{2}(\cdots )+\cdots \over 
\alpha _{1}^{2}+\cdots +\alpha _{p}^{2}} \le \lambda _{1}^{2}$$
and equality holds if and only if $\alpha _{r+1}=\cdots =\alpha _{p}=0$,
i.e. if and only if ${\bf x}$ belongs to the eigenspace corresponding 
to $\lambda _{1}$. \quad {\hspace* {\fill}{$\Box $}}
\medskip 

\noindent 
Note that analogous statement like Theorem 2.3 holds also if we consider 
${\bf x}$ as vector of $\Sigma _{2}$ and $\varphi $ is the angle between 
${\bf x}$ and $\Sigma _{1}$. Thus we obtain the following geometrical 
interpretation: 

Among all values $\cos^{2} \varphi $ 
where $\varphi $ is angle between any vector ${\bf x}\in \Sigma _{1}$ 
and any vector ${\bf y}\in \Sigma _{2}$, the maximal value $\lambda _{1}^2$ 
is the first (maximal) principal value. Then 
$$\Sigma _{11}=\{ {\bf x}\in \Sigma _{1}\vert \cos ^{2}
({\bf x},\Sigma _{2}) = \lambda _{1}^{2}\}$$
is the the principal subspace of $\Sigma _{1}$. Analogously 
$$\Sigma _{21}=\{ {\bf y}\in \Sigma _{2}\vert 
\cos ^{2}({\bf y},\Sigma _{1}) = \lambda _{1}^{2}\}$$
is the principal subspace of $\Sigma _{2}$ and moreover 
$dim \Sigma _{11}=dim \Sigma _{21}$. 
Now let us consider the subspaces $\Sigma '_1$ and $\Sigma '_2$ where 
$\Sigma '_1$ is orthogonal complement of $\Sigma _{11}$ in 
$\Sigma _{1}$ and $\Sigma '_2$ is orthogonal complement of 
$\Sigma _{21}$ in $\Sigma _{2}$. 
Among all values $\cos^{2} \varphi $ 
where $\varphi $ is angle between any vector ${\bf x}\in \Sigma '_{1}$ 
and any vector ${\bf y}\in \Sigma '_{2}$, the maximal value $\lambda _{2}^2$ 
is the second principal value. Then 
$$\Sigma _{12}=\{ {\bf x}\in \Sigma '_{1}\vert 
\cos ^{2}({\bf x},\Sigma '_{2}) = \lambda _{2}^{2}\}$$
is the principal subspace of $\Sigma '_{1}$. Analogously 
$$\Sigma _{22}=\{ {\bf y}\in \Sigma '_{2}\vert 
\cos ^{2}({\bf y},\Sigma '_{1}) = \lambda _{2}^{2}\}$$
is the principal subspace of $\Sigma '_{2}$ and moreover 
$dim \Sigma _{12}=dim \Sigma _{22}$. Continuing this procedure we obtain 
the decompositions of orthogonal principal subspaces 
$$\Sigma _1 = \Sigma _{11} + \Sigma _{12} + \cdots + \Sigma _{1,s+1} $$
$$\Sigma _2 = \Sigma _{21} + \Sigma _{22} + \cdots + \Sigma _{2,s+1} $$
where $dim \Sigma _{1i}=dim \Sigma _{2i}$, $(1\le i\le s)$. 
The subspaces $\Sigma _{1,s+1}$ and $\Sigma _{2,s+1}$ correspond for the 
possible value 0 as a principal value. 

\medskip
\noindent
{\bf Example.} Let $\Sigma _{1}$ be generated by the vectors 
$(1,0,0,0)$ and $(0,1,0,0)$ and $\Sigma _{2}$ be generated by 
$(\cos \varphi ,0,\sin \varphi ,0)$ and $(0,\cos \varphi ,0,\sin \varphi )$. 
Then $\cos ^{2}\varphi $ is unique principal value, its multiplicity is 2 
and $\Sigma _1$ and $\Sigma _{2}$ are principal subspaces themselves. 

\medskip
\noindent
At the end we prove a theorem which determines the 
orthogonal projection of any vector ${\bf x}$ on any subspace of $E_{n}$. 

\medskip 
\noindent 
{\bf Theorem 2.4.} {\em In the $n$-dimensional Euclidean space $E_{n}$ 
let be given a subspace $\Sigma $ generated by 
$k$ linearly independent vectors ${\bf a}_{i}$, 
$(1\le i\le k;\; k\le n-1)$. The orthogonal projection ${\bf x}'$ of 
an arbitrary vector ${\bf x}$ of $E_{n}$ is given by }
$$ {\bf x}' = -{1\over \Gamma } \left \vert \matrix{
{\bf 0} & ({\bf x},{\bf a}_{1}) & ({\bf x},{\bf a}_{2}) & \cdots &
({\bf x},{\bf a}_{k})\cr 
{\bf a}_{1} & ({\bf a}_{1},{\bf a}_{1})&({\bf a}_{1},{\bf a}_{2})&\cdots &
({\bf a}_{1},{\bf a}_{k})\cr 
{\bf a}_{2} & ({\bf a}_{2},{\bf a}_{1})&({\bf a}_{2},{\bf a}_{2})&\cdots &
({\bf a}_{2},{\bf a}_{k})\cr 
\cdot \cr 
\cdot \cr 
\cdot \cr 
{\bf a}_{k} & ({\bf a}_{k},{\bf a}_{1})&({\bf a}_{k},{\bf a}_{2})&\cdots &
({\bf a}_{k},{\bf a}_{k})\cr } \right \vert , \leqno{(2.1)}$$
{\em where $\Gamma $ is the Gram's determinant of the vectors }
${\bf a}_{i}$, $(1\le i\le k)$.

\medskip
\noindent {\em Proof.} 
According to (2.1) it is obvious that 
$$ {\bf x}-{\bf x}' = {1\over \Gamma } \left \vert \matrix{
{\bf x} & ({\bf x},{\bf a}_{1}) & ({\bf x},{\bf a}_{2}) & \cdots &
({\bf x},{\bf a}_{k})\cr 
{\bf a}_{1} & ({\bf a}_{1},{\bf a}_{1})&({\bf a}_{1},{\bf a}_{2})&\cdots &
({\bf a}_{1},{\bf a}_{k})\cr 
{\bf a}_{2} & ({\bf a}_{2},{\bf a}_{1})&({\bf a}_{2},{\bf a}_{2})&\cdots &
({\bf a}_{2},{\bf a}_{k})\cr 
\cdot \cr 
\cdot \cr 
\cdot \cr 
{\bf a}_{k} & ({\bf a}_{k},{\bf a}_{1})&({\bf a}_{k},{\bf a}_{2})&\cdots &
({\bf a}_{k},{\bf a}_{k})\cr }\right \vert .$$

By scalar multiplication of this equality by ${\bf a}_{i}$, $(1\le i\le k)$
the first column is equal to the $(i+1)$-st column and thus 
$$ ({\bf x}-{\bf x}',{\bf a}_{i}) = 0, \quad (1\le i\le k).$$
Since ${\bf x}'$ is linear combination of the vectors ${\bf a}_{i}$, 
$(1\le i\le k)$ then the vector ${\bf x}'$ lies in $\Sigma $. Moreover, 
${\bf x}-{\bf x}'$ is orthogonal to the base vectors of $\Sigma $, we 
obtain that ${\bf x}'$ is the required orthogonal projection of 
${\bf x}$ on the subspace ${\Sigma }$. {\hspace* {\fill}{$\Box $}}

\medskip
\noindent {\bf 3. Principle of duality and canonical form}
\medskip

In this section we will consider the duality principle like 
in the Theorem 1.3
and as a crown of all previous research will be given the canonical 
form of two subspaces $\Sigma _{1}$ and $\Sigma _{2}$. 
Now let $\Sigma _{i}^{*}$ denote the orthogonal subspace 
of $\Sigma _{i}$, $(i=1,2)$ in the Euclidean space $E_{n}$. We saw that 
$\varphi (\Sigma _{1},\Sigma _{2})=\varphi (\Sigma _{1}^* ,\Sigma _{2}^* )$
and now the same conclusions for the eigenvalues and principal 
subspaces (principal directions) also hold for the subspaces 
$\Sigma _{1}^* $ and $\Sigma _{2}^* $. 

\medskip
\noindent
{\bf Theorem 3.1.} {\em If $\Sigma _1 $ and $\Sigma _2 $ are any subspaces
of the Euclidean vector space $E_{n}$ and 
$\Sigma _1 ^* $ and $\Sigma _2 ^*$ are their orthogonal complements, 
then the nonzero and different from 1 
principal values for the pair $(\Sigma _{1},\Sigma _{2})$ are the same 
for the pair $(\Sigma _{1}^{*},\Sigma _{2}^{*})$ with the same 
multiplicities and conversely. 

If $p + q \le n$, then the multiplicity of 1 for the pair 
$(\Sigma _{1}^{*},\Sigma _{2}^{*})$ is bigger for $n-p-q$ 
than the multiplicity of 1 for the pair $(\Sigma _{1},\Sigma _{2})$.

If $p + q \ge n$, then the multiplicity of 1 for the pair 
$(\Sigma _{1},\Sigma _{2})$ is bigger for $p+q-n$ than the
multiplicity of 1 for the pair $(\Sigma _{1}^{*},\Sigma _{2}^{*})$. }

\medskip
\noindent {\em Proof.} 
We use the same notations and assumptions as in the proof of 
the Theorem 1.3. Specially, the matrices $X$, $Y$ and $Z$ are the same. 
Assume that $p + q \le n$. The case $n>p + q $ can be 
discussed analogously. 

We will prove the following identity 
$$det [\lambda I_{p\times p} - XX^{T}]\cdot 
(\lambda -1)^{n-q-p}= 
det [\lambda I_{(n-q)\times (n-q)} - Z^{T}Z]$$
and hence the proof will be finished. 

Since $A$ is an orthogonal matrix, it holds 
$$XX^{T} = I_{p \times p } - YY^{T}\quad \hbox { and } \quad 
Z^{T}Z = I_{(n-q) \times (n-q )} - Y^{T}Y$$
and we should prove that 
$$det [(\lambda -1)I_{p \times p } + YY^{T}]\cdot 
(\lambda -1)^{n-q-p} = 
det [(\lambda -1)I_{(n-q) \times (n-q)} + Y^{T}Y].$$
Multiplying this equality by $(-1)^{n-q}$ and putting $1-\lambda =\mu $, 
we should prove that 
$$det [\mu I_{p \times p } - YY^{T}] \cdot \mu ^{n-q-p}= 
det [\mu I_{(n-q) \times (n-q)} - Y^{T}Y].$$
Let $\mu _{1},\cdots ,\mu _{p}$ be the eigenvalues of $YY^{T}$. 
According to Theorem 1.2, both sides of the last equality are 
equal to 
$$(\mu  -\mu_{1})(\mu  -\mu  _{2})\cdots  (\mu  -\mu  _{p})\mu 
^{n-q-p }.\eqno{\Box}$$

\noindent According to Theorem 3.1 we obtain the following consequence. 

\medskip
\noindent
{\bf Corollary 3.2.} {\em According to the notations of the Theorem 3.1, 

i) the number of nonzero and nonunit principal values 
(each value counts as many times as its multiplicity) of the pair 
$(\Sigma _1,\Sigma _2)$ is less or equal to $n/2$; 

ii) if $n$ is odd number and $p=q$, then at least one of the pairs 
$(\Sigma _1,\Sigma _2)$ and $(\Sigma _1^*,\Sigma _2^*)$ has 
a principal value 1, i.e. they have a common subspace of 
dimension $\ge 1$.} 

\medskip
\noindent
Now we are able to give the canonical form of two subspaces. 
In order to avoid many indices we assume that the 
considered subspaces of $E_{n}$ are $\Sigma $ and $\Pi$ with dimensions 
$p$ and $q$ respectively. We denote by $\Sigma ^*$ and $\Pi ^*$ 
the orthogonal subspaces of $E_{n}$. Without loss of generality 
we assume that $p\le q$. Since the canonical form is according to these 
four subspaces, we can also assume that $p+q\le n$. Indeed, if $p+q>n$ 
then $(n-p)+(n-q)<n$ and 
we can consider the subspaces $\Sigma ^*$ and $\Pi ^*$. 

Assume that $1=c_{0}>c_{1}>c_{2}>\cdots >c_{s}>c_{s+1}=0$ be the 
principal values for the pair $(\Sigma ,\Pi )$ with multiplicities 
$r_{0},r_{1},\cdots ,r_{s+1}$ respectively, such that 
$p=r_{0}+r_{1}+\cdots +r_{s+1}$. Let $\Sigma $ be generated by the 
following orthonormal vectors 
$${\bf a}_{01},\cdots ,{\bf a}_{0r_{0}}, 
{\bf a}_{11},\cdots ,{\bf a}_{1r_{1}}, \cdots ,
{\bf a}_{s1},\cdots ,{\bf a}_{sr_{s}}, 
{\bf a}_{s+1,1},\cdots ,{\bf a}_{s+1,r_{s+1}},$$ 
such that the vectors ${\bf a}_{i1},\cdots ,{\bf a}_{ir_{i}}$ 
generate the principal subspace for the principal value $c_{i}$, 
$(0\le i\le s+1)$. The pair of subspaces $(\Sigma ^*,\Pi ^*)$ have 
the same principal values $1=c_{0}>c_{1}>c_{2}>\cdots >c_{s}>c_{s+1}=0$
with multiplicities 
$r'_{0}=r_{0}+n-p-q,r_{1},\cdots ,r_{s+1}$. 
Assume that $\Sigma ^*$ 
is generated by the following orthonormal vectors 
$${\bf a}^{*}_{01},\cdots ,{\bf a}^{*}_{0r'_{0}}, 
{\bf a}^{*}_{11},\cdots ,{\bf a}^{*}_{1r_{1}}, \cdots ,
{\bf a}^{*}_{s1},\cdots ,{\bf a}^{*}_{sr_{s}}, 
{\bf a}^{*}_{s+1,1},\cdots ,{\bf a}^{*}_{s+1,r_{s+1}},
{\bf a}^{*}_{1},\cdots ,{\bf a}^{*}_{q-p}$$
where the vectors ${\bf a}_{i1},\cdots ,{\bf a}_{ir_{i}}$ 
generate the principal subspace for the principal value $c_{i}$, 
$(1\le i\le s+1)$, ${\bf a}_{01},\cdots ,{\bf a}_{0r'_{0}}$ 
generate the principal subspace for the principal value 1 and 
${\bf a}^{*}_{1},\cdots ,{\bf a}^{*}_{q-p}$ 
be the remaining $q-p$ orthonormal vectors. 

Now we choose the orthonormal vectors of $\Pi $ as follows. We choose 
$${\bf b}_{01},\cdots ,{\bf b}_{0r_{0}}, 
{\bf b}_{11},\cdots ,{\bf b}_{1r_{1}}, \cdots ,
{\bf b}_{s1},\cdots ,{\bf b}_{sr_{s}}, 
{\bf b}_{s+1,1},\cdots ,{\bf b}_{s+1,r_{s+1}},
{\bf b}_{1},\cdots ,{\bf b}_{q-p} $$ 
such that ${\bf b}_{0i}$ coincides with ${\bf a}_{0i}$, $(1\le i\le r_{0})$,
${\bf b}_{i1},\cdots ,{\bf b}_{ir_{i}}$ 
generate the principal subspace for the principal value $c_{i}$, 
$(1\le i\le s)$ and such that 
$({\bf a}_{iu},{\bf b}_{iv})=\delta _{uv}c_{i}$. 
The vectors ${\bf b}_{s+1,1},\cdots ,{\bf b}_{s+1,r_{s+1}}$ generate the 
same subspace as the vectors 
${\bf a}^*_{s+1,1},\cdots ,{\bf a}^*_{s+1,r_{s+1}}$ and we can choose 
${\bf b}_{s+1,i}={\bf a}^*_{s+1,i}$, $(1\le i\le r_{s+1})$. The 
vectors ${\bf b}_{1},\cdots ,{\bf b}_{q-p}$ generate the same 
space as the vectors ${\bf a}^*_{1},\cdots ,{\bf a}^*_{q-p}$ 
and we can choose ${\bf b}_{i}={\bf a}^*_{q-p+1-i}$, $(1\le i\le q-p)$. 

Finally we determine the orthonormal vectors of $\Pi ^*$ 
$${\bf b}^{*}_{01},\cdots ,{\bf b}^{*}_{0r'_{0}}, 
{\bf b}^{*}_{11},\cdots ,{\bf b}^{*}_{1r_{1}}, \cdots ,
{\bf b}^{*}_{s1},\cdots ,{\bf b}^{*}_{sr_{s}}, 
{\bf b}^{*}_{s+1,1},\cdots ,{\bf b}^{*}_{s+1,r_{s+1}}$$
as follows. The vectors ${\bf b}^{*}_{01},\cdots ,{\bf b}^{*}_{0r'_{0}}$ 
can be chosen such that ${\bf b}^{*}_{0i}={\bf a}^{*}_{0i}$, 
$(1\le i\le r'_{0})$. 
The vectors ${\bf b}^*_{i1},\cdots ,{\bf b}^*_{ir_{i}}$ 
generate the principal subspace for the principal value $c_{i}$, 
$(1\le i\le s)$, and the vectors 
${\bf b}^*_{i1},\cdots ,{\bf b}^*_{ir_{i}}$ 
can uniquely be chosen such that 
$({\bf a}^*_{iu},{\bf b}^*_{iv})=\delta _{uv}c_{i}$. The vectors 
${\bf b}^{*}_{s+1,1},\cdots ,{\bf b}^{*}_{s+1,r_{s+1}}$ 
generate the same subspace as the vectors 
${\bf a}^{*}_{s+1,1},\cdots ,{\bf a}^{*}_{s+1,r_{s+1}}$ and thus we can 
choose ${\bf b}^{*}_{s+1,i}={\bf a}^{*}_{s+1,i}$, $(1\le i\le r_{s+1})$. 

Moreover, the vectors 
${\bf a}^{*}_{11},\cdots ,{\bf a}^{*}_{1r_{1}}, \cdots ,
{\bf a}^{*}_{s1},\cdots ,{\bf a}^{*}_{sr_{s}}$ 
can be chosen such that 
$$({\bf a}^*_{iu},{\bf b}_{iv})=-\delta _{uv}\sqrt{1-c_{i}^{2}}, 
\quad (1\le i\le s).$$
Now we know some of the inner products 
between the base vectors of $\Sigma $ and $\Sigma ^*$ and the base vectors 
of $\Pi $ and $\Pi ^*$. The matrix $P$ of all such $n\times n$ inner 
products must be orthogonal and can uniquely be obtained from the 
above inner products. Considering the base vectors of 
$\Sigma $ in the mentioned order together with the base vectors of 
$\Sigma ^*$ in the opposite order 
and on the other side the base vectors of $\Pi $ in the mentioned order 
together with the base vectors of $\Pi ^*$ in the opposite order 
we obtain the following 
$$(r_{0}+r_{1}+r_{2}+\cdots +r_{s}+r_{s+1}+(q-p)+
r_{s+1}+r_{s}+\cdots +r_{2}+r_{1}+r'_{0})\times $$
$$\times (r_{0}+r_{1}+r_{2}+\cdots +r_{s}+r_{s+1}+(q-p)+
r_{s+1}+r_{s}+\cdots +r_{2}+r_{1}+r'_{0})$$ 
matrix as {\em canonical matrix} for the subspaces $\Sigma $ and $\Pi $: 
$$ P = \left [\matrix{
I&0&0&\cdots &0&0&0&0&0&\cdots &0&0&0\cr 
0&c_{1}I&0&\cdots &0&0&0&0&0&\cdots &0&
d_{1}I'&0\cr 
0&0&c_{2}I&\cdots &0&0&0&0&0&\cdots &
d_{2}I'&0&0\cr 
\cdot \cr
\cdot \cr
\cdot \cr
0&0&0&\cdots &c_{s}I&0&0&0&
d_{s}I'&\cdots &0&0&0\cr 
0&0&0&\cdots &0&0&0&I'&0&\cdots &0&0&0\cr 
0&0&0&\cdots &0&0&I&0&0&\cdots &0&0&0\cr 
0&0&0&\cdots &0&I'&0&0&0&\cdots &0&0&0\cr 
0&0&0&\cdots &-d_{s}I'&0&0&0&
c_{s}I&\cdots &0&0&0\cr 
\cdot \cr 
\cdot \cr
\cdot \cr 
0&0&-d_{2}I'&\cdots &0&0&0&0&0&
\cdots &c_{2}I&0&0\cr 
0&-d_{1}I'&0&\cdots &0&0&0&0&0&
\cdots &0&c_{1}I&0\cr 
0&0&0&\cdots &0&0&0&0&0&
\cdots &0&0&I\cr }\right ] ,$$
where $d_{i}=\sqrt{1-c_{i}^{2}}$, $(1\le i\le s)$ and 
$I'$ denotes the matrix with 1 on the opposite diagonal of the main
diagonal and the other elements are zero. 

Note that the principal values for the pair $(\Sigma ,\Pi ^*)$ 
(also $(\Sigma ^*,\Pi ))$ are the numbers 
$d_{i}^{2}=1-c_{i}^{2}=\sin ^{2}\varphi _{i}$ with the same multiplicities 
as $c_{i}^{2}$. Moreover the previous canonical matrix $P$ is also 
canonical matrix for the pair $(\Sigma ,\Pi ^*)$ (also $(\Sigma ^*,\Pi ))$ 
if we permute its rows and columns. Then the order 
$q-p$ converts into $n-p-q$ and vice versa. 

The previous consideration yields to the following statement. 

\medskip
\noindent
{\bf Theorem 3.3.} {\em 
Let $n,p,q$ be positive integers such that $n\le p+q$ and $p\le q$. 
Then for any $p$ values $c_{1}^{2},\cdots ,c_{p}^{2}$, 
$(0\le c_{i} \le 1)$ 
there exist two subspaces $\Sigma _{1}$ and $\Sigma _{2}$ of $E_n$
with dimensions $p$ and $q$ such that $c_{1}^{2},\cdots ,c_{p}^{2}$ are 
principal values for the pair $(\Sigma _{1},\Sigma _{2})$. 
The existence of the 
subspaces $\Sigma _1$ and $\Sigma _2$ is uniquely 
up to orthogonal motion in $E_n$. }

\medskip
\noindent {\em Proof.} 
Let $n,p,q$ be positive integers such that $n\le p+q$ and $p\le q$ and 
let be given $p$ values $c_{i}^{2}$, $(0\le c_{i} \le 1)$. 
We choose arbitrary orthonormal base 
${\bf a}_{1},\cdots ,{\bf a}_{p},{\bf a}^{*}_{n-p},\cdots ,{\bf a}^{*}_{1}$
of $E_n$. Then we introduce $q$ vectors ${\bf b}_{1},\cdots ,
{\bf b}_{q}$ whose coordinates with respect to 
${\bf a}_{1},\cdots ,{\bf a}_{p},{\bf a}^{*}_{n-p}$, 
$\cdots ,{\bf a}^{*}_{1}$
are given by the first $q$ 
columns of the matrix $P$. Then it is obvious that the principal values 
for the pair $(\Sigma _{1},\Sigma _{2})$ where $\Sigma _{1}$ is generated by 
${\bf a}_{1},\cdots ,{\bf a}_{p}$ and $\Sigma _{2}$ is generated by the 
vectors ${\bf b}_{1},\cdots ,{\bf b}_{q}$ are just the given numbers 
$c_{1}^{2},\cdots ,c_{p}^{2}$. 

Let $(\Sigma _1,\Sigma _2)$ and $(\Sigma '_1,\Sigma '_2)$ be two pairs of 
subspaces with the same principal values. Without loss of generality 
we assume that both of them are given in canonical form given by the 
same canonical matrix $P$. Let 
$$\{{\bf a}_{1},\cdots ,{\bf a}_{p},{\bf a}^*_{1},\cdots ,{\bf a}^*_{n-p}\}
\quad \hbox { and }\quad 
\{{\bf a}'_{1},\cdots ,{\bf a}'_{p},
{\bf a}'^*_{1},\cdots ,{\bf a}'^*_{n-p}\}$$
be the base vectors of $\Sigma _{1}+\Sigma ^*_{1}$ and 
$\Sigma '_{1}+\Sigma '^*_{1}$ corresponding to their canonical forms. 
Since the base vectors of $\Sigma _{2}+\Sigma ^*_{2}$ and 
$\Sigma '_{2}+\Sigma '^*_{2}$ are determined uniquely, it is sufficient to
choose the orthogonal transformation 
$\varphi $ which maps the mentioned base of 
$\Sigma _{1}+\Sigma ^*_{1}$ into the mentioned base of 
$\Sigma '_{1}+\Sigma '^*_{1}$ and then 
$\varphi (\Sigma _{1})=\Sigma '_{1}$ and 
$\varphi (\Sigma _{2})=\Sigma '_{2}$. \quad {\hspace* {\fill}{$\Box $}}

\medskip
\noindent
{\bf Theorem 3.4.} {\em Let $A$ be a symmetric matrix of $n$-th order. 
Assume that the linear subspace $L$ of $E_{n}$ such that $A$ is positive 
definite matrix in $L$ and $A^{-1}$ is positive definite matrix in 
the orthogonal complement $L^{*}$, then $A$ is positive definite matrix. }

\medskip
\noindent {\em Proof.} 
If $A\vert L$ denotes the restriction of $A$ to $L$, and 
by $ind (A\vert L)$ is denoted the number of negative eigenvalues of 
$V^{T}AV$, where $V$ is the matrix of the base of $L$, then the following 
lemma holds. 

\medskip
\noindent
{\bf Lemma 3.5.} {\em Let $A$ be a symmetric nonsingular matrix of $n$-th 
order, 
and let $L$ and $L^{*}$ be the same notations like in Theorem 3.4. 
If $A^{-1}\vert L^{*}$ is nonsingular restriction, then also the 
restriction $A\vert L$ is nonsingular and moreover }
$$ ind(A\vert E_{n}) = ind (A\vert L) + ind (A^{-1}\vert L^{*}).$$

\noindent The Theorem 3.4 obtains for the special case 
$$ ind(A\vert L)=ind(A^{-1}\vert L^{*}) = 0.$$

\noindent 
{\em Proof of the Lemma 3.5.} Let $V$ and $W$ denote the matrices from the 
bases of $L$ and $L^{*}$ respectively. Then $B=AVW$ is 
nonsingular matrix. Indeed, it is supposed that $AV{\bf x}=W{\bf y}$ 
for the vectors ${\bf x}$ and ${\bf y}$. Multiplying this equality 
by $W^{*}A^{-1}$ from left, we obtain $W^{*}A^{-1}W{\bf y}=0$, because 
$V^{*}W=0$. This implies ${\bf y}=0$ which means that 
$W^{*}A^{-1}W$ is nonsingular matrix. Consequently, 
$W{\bf x}=A^{-1}W{\bf y}=0$ implies ${\bf x}=0$. It implies that 
$$ind (A\vert E_{n}) = ind (A^{-1}\vert E_{n}) = ind (B^{T}A^{-1}B\vert E_n)
=$$
$$= ind (A\vert L) + ind (A^{-1}\vert L^{*}) .
\eqno{\Box}$$

\end{document}